\documentclass[12pt,a4paper]{article}
\usepackage{mathrsfs}

\usepackage{amssymb, latexsym, amsthm}
\usepackage{graphicx} % avoid epsfig or earlier such packages
\usepackage{url}      % for URLs and DOIs
\usepackage{amsmath}  % many want amsmath extensions
\IfFileExists{ajr.sty}{\usepackage{ajr}}{}

\newcommand{\Ord}[1]{\ensuremath{\mathcal O\big(#1\big)}}
\newcommand{\rat}[2]{{\textstyle\frac{#1}{#2}}}
\newcommand{\doi}[1]{\href{http://dx.doi.org/#1}{doi:#1}}
\newcommand{\Z}[1]{\ensuremath{e^{#1 t}\star}}
\newcommand{\Za}[1]{\ifcase#1 undef0
  \or undef1
  \or\Z{-\rat{27}{10}}
  \or\Z{-\rat{38}{5}}
  \else undefinf \fi}
\newcommand{\Zb}[1]{\ifcase#1 undef0
  \or\Z{-\rat{2}{\epsilon}}
  \or\Z{-\rat{5}{\epsilon}}
  \or\Z{-\rat{10}{\epsilon}}
  \else undefinf \fi}

\newcommand{\p}{\partial}

\renewcommand{\phi}{\varphi}
\newcommand{\e}{\epsilon}

\newcommand{\R}{{\mathbb R}}

\newcommand{\PX}{{\Bbb{P}}}

\newcommand{\ssm}{\textsc{ssm}}
\newcommand{\spde}{\textsc{spde}}
\newcommand{\sde}{\textsc{sde}}
\newcommand{\ldp}{\textsc{ldp}}

\newtheorem{theorem}{Theorem}
\newtheorem{lemma}[theorem]{Lemma}
\newtheorem{corollary}[theorem]{Corollary}
\theoremstyle{definition}
\newtheorem{definition}[theorem]{Definition}
\newtheorem{remark}[theorem]{Remark}

\title{Large deviations for slow-fast stochastic partial differential equations}

\author{
Wei Wang\thanks{School of Mathematical Sciences, University of
Adelaide, Adelaide, \textsc{Australia}. \protect\url{mailto:
w.wang@adelaide.edu.au}; and Department of Mathematics, Nanjing
University, Nanjing, \textsc{China}.
\protect\url{mailto:wangweinju@yahoo.com.cn} } \and A.~J.
Roberts\thanks{School of Mathematical Sciences, University of
Adelaide, Adelaide, \textsc{Australia}.
\protect\url{mailto:anthony.roberts@adelaide.edu.au}} \and Jinqiao
Duan\thanks{Department of Applied Mathematics, Illinois Institute of
Technology, Chicago, IL 60616, \textsc{usa}
\protect\url{mailto:duan@iit.edu}}
 }

\date{\today}

\begin{document}

% Use default \verb|\maketitle|.
\maketitle

% Use the \verb|abstract| environment.
\begin{abstract}
A large deviation principle is derived for stochastic partial
differential equations with slow-fast components. The result shows
that the rate function is exactly that of the averaged equation plus
the fluctuating deviation which is a stochastic partial differential
equation with small Gaussian perturbation. This also confirms the
effectiveness of the approximation of the averaged equation plus the
fluctuating deviation to the slow-fast stochastic partial
differential equations.
\end{abstract}

%\tableofcontents

\section{Introduction}\label{sec:intro}

Uncertainties (noise) is widely recognized in modeling, analyzing,
simulating and predicting complex phenomena~\cite[e.g.]{Arnold01,
E00, Imkeller, TM05}. Noise causes rare events in nonlinear
stochastic system describing the metastability of the
system~\cite{ERE, KN85, FW98}. The classic example is a tunnelling
event between two stable points of a macroscopic system. This rare
event eventually occurs after a long time scale by an addition of a
small external noise to the macroscopic system, but the probability
of such rare event converges to zero as the strength of noise tends
to zero. We need to understand the rate of such convergence. The
theory of large deviations provide a powerful tool to give an
estimate to the rate of such convergence, shown to be exponential
for finite dimensional stochastic systems~\cite[e.g.]{Var, FW98,
Stroock}. To study the metastability of a macroscopic system with
small noise we must build the large deviation principle~(\ldp).

Stochastic partial differential equations (\spde s) are appropriate mathematical models for many multiscale systems with uncertain influences~\cite{WaymireDuan}.
The \ldp\ for \spde s  has been studied by many people~\cite[e.g.]{BDM08,CeRoc,CM97, Chow, Duan09, KX96, Pes, PZ92, Sow92} under a different framework to that used here.
However, there are very few results on the \ldp\ for nonlinear stochastic system with two widely separated timescales, as often appears in a complex system.
Freidlin and Wentzell~\cite{FW98} first studied the \ldp\ for finite dimensional slow-fast stochastic system with partial coupling.
They used a bounded assumption on the nonlinearity.
Then Veretennikov~\cite{Yu99} built a \ldp\ for the full coupled case with bounded assumptions on nonlinearity.
Later still, under a bounded assumption, Ioffe studied the \ldp\ for stochastic reaction-diffusion equation with rapidly oscillating random noise in the special case when there is no coupling between the  slow component  and fast component~\cite{Iof, Iof91}.
However, there appear to be no other previous \ldp{}s for slow-fast coupled stochastic partial differential equations.

This article establishes the Freidlin--Wentzell \ldp{} for a class of \spde s with stochastic fast component and deterministic slow component.
Let~$D$ be an open bounded interval and $L^2(D)$ be the Lebesgue space of square integrable real valued functions on~$D$.
Consider the following pair of stochastically forced, coupled, reaction-diffusion \spde{}s for any~$\e>0$
 \begin{eqnarray}
&&\p_t u^\e=\p_{xx}u^\e+f(u^\e, v^\e)\,,\quad u^\e(0)=u_0\in L^2(D) \label{e:ue}\\
&&\p_tv^\e=\frac{1}{\e}\big[\p_{xx}v^\e+g(u^\e, v^\e)]+
   \frac{\sigma}{\sqrt{\epsilon}}\p_tW(t)\,, \quad v^\e(0)=v_0\in L^2(D)
   \label{e:ve}
\end{eqnarray}
with zero Dirichlet boundary condition on~$\partial D$.
Here $W(t)$~is an  $L^2(D)$~valued Wiener process defined on a complete probability space $(\Omega, \mathcal{F}, \mathbb{P})$ as detailed in the next section.
In physical applications this supposition is that the noise directly drives microscopic modes~$v^\epsilon$---the noise only emerges in the macroscopic modes~$u^\epsilon$ through nonlinear coupling.

If parameter $\e>0$ is small so that $v^\e$~fluctuates rapidly, then an effective approximated system is desirable.
Under some appropriate assumptions one averages the slow part~$u^\e$ over the fast part~$v^\e$ which yields the following so-called averaged equation describing the dynamics of the system on a slow time scale
\begin{equation}\label{e:averaged}
\p_tu=\p_{xx}u+\bar{f}(u)\,, \quad u(0)=u_0\,,\quad u|_{\p D}=0\,.
\end{equation}
Here $\bar{f}(u)$~is the average of~$f(u,v)$ over the distribution of the fast part~$v$.
Cerrai et al.~\cite[e.g.]{CF09,C09} recently developed more on  the averaging principle of \spde s.
Wang and Roberts~\cite{WR08} very recently gave a further approximation result via a martingale discussion  which shows that the deviation $u^\e(t)-u(t)$ is approximated by~$\sqrt{\e}z(t)$, for some Guassian process~$z(t)$, in the case that fast component is coupled with slow component and without any Lipschitz assumption on the slow component.
Cerria~\cite{Cer09} obtained the same result for the special case where the nonlinearity is Lipschitz and there is no coupling of the slow component to the fast component.
Then the deviation estimate shows that~$u^\e$ approximates~$u$ with a small Gaussian perturbation and this suggests a \ldp{} for~$\{u^\e\}_\e$.
By studying the \ldp{} for some auxiliary systems, we prove the \ldp{} for~$\{u^\e\}_\e$, Theorem \ref{thm:main}.
Moreover, the rate function for the \ldp{} of $\{u^\e\}_\e$ in the main result is exactly that of~$\{\tilde{u}^\e\}_\e$ solving~(\ref{e:au1+deviation})--(\ref{e:au2+deviation}), which is the averaged equation~(\ref{e:averaged}) plus deviation up to errors of~$\mathcal{O}(\e)$, Section~\ref{sec:pre}.
Our results further confirms the effectiveness of the averaged equation plus deviation to approximate slow-fast \spde{}s~(\ref{e:ue})--(\ref{e:ve}).

Recently, a weak convergence approach, which avoids giving some technical exponential tightness estimates, was applied to obtain \ldp\ for \spde s~\cite[e.g.]{DE97, WD09}.
But this approach does not work here because a drift transformation leads the fast system to become a non-autonomous system for which one cannot average the slow part over the fast part.
For this here we still give some exponential tight estimates, Section~\ref{sec:exp tight}, and then by some contraction principles and an approximation we obtain \ldp\ for some auxiliary systems, Section~\ref{sec:aux-sys}.
An approximation shows that the slow-fast stochastic system~(\ref{e:ue})--(\ref{e:ve}) is comparable with the auxiliary systems near some functions, Section~\ref{sec:ldp}, and we derive the \ldp{} for~$\{u^\e\}_\e$.

Section~\ref{sec:esfsrd} presents an example slow-fast reaction-diffusion \spde\ to illustrate the \ldp\ theory.
Section~\ref{sec:ssmeg} then explores the parameter regime near a stochastic pitchfork bifurcation in this example.
Constructing the stochastic `superslow' manifold and the evolution thereon confirms that there is indeed a close correspondence between the original example system and the \ldp\ averaged system.

%One key step here is to prove \ldp{} for some auxiliary slow-fast
%stochastic systems, Section~\ref{sec:aux-sys}\,. Then by an
%approximating discussion which shows that the slow-fast stochastic
%system~(\ref{e:ue})--(\ref{e:ve}) is comparable with the auxiliary
%system near any function, Section~\ref{sec:ldp}\,.
%
%
%Our result shows that rate function of the slow-fast \spde{} is
%exactly that of the averaged equation plus the small fluctuation
%(\ref{e:au1+deviation})--(\ref{e:au2+deviation})\,.

\section{Preliminaries}\label{sec:pre}

Let $H=L^2(D)$ with $L^2$-norm denoted by~$\|\cdot\|_0$ and inner product by~$\langle\cdot, \cdot\rangle$.
Define the linear operator $A=\p_{xx}$ with zero Dirichlet boundary condition on~$D$.
Then operator~$A$ is the generator of a compact analytic semigroup~$e^{At}$, $t\geq 0$\,, on~$H$.
Moreover, denote by~$\{e_i\}_{i=1}^\infty$\,, which forms a complete standard orthogonal basis in~$H$, a family of eigenfunctions of~$A$  and $-Ae_i=\lambda_ie_i$\,, $\lambda_i>0$\,, $i=1, 2, \ldots$\,.
For any $\alpha>0$ and $u\in H$ define $\|u\|_\alpha=\|A^{\alpha/2}u\|_0$\,.
Then let $H_0^\alpha$~be the space of the closure of~$C_0^\infty(D)$, the space of smooth functions with compact support on~$D$, under the norm~$\|\cdot\|_\alpha$.
Furthermore, let $H^{-\alpha}$~denote the dual space of~$H_0^\alpha$\,.
Also we are given $H$~valued Wiener processes~$W(t)$, $t\geq 0$\,, defined on the complete probability space~$(\Omega, \mathcal{F}, \mathcal{F}_t, \PX)$~\cite{PZ92}.
Denote by~$\mathbb{E}$ the expectation operator with respect to~$\PX$.
We consider the \spde{}s of the form (\ref{e:ue})--(\ref{e:ve})  with separated time scale and with $\sigma\neq 0$ is an arbitrary real number parametrising the strength of the noise.  We adopt the following four hypotheses.

\begin{description}
  \item[H$_1$] $f(x, y):\R\times \R\rightarrow\R$ is Lipschitz continuous in both $x$ and $y$ with Lipschtiz constant~$L_f$ and for all $x,y\in\R$
 \begin{equation*}
 |f(x, y)|^2\leq ax^2+by^2+c \,,
 \quad  f(x, y)x\leq a x^2+bxy+c \,,
 \end{equation*}
  for some positive constants~$a$, $b$ and~$c$.
 \item[H$_2$] $g(x, y):\R\times\R\rightarrow\R$ is Lipschitz continuous in both variables with Lipschitz constant~$L_g$.
For any $x,y\in\R$
 \begin{equation*}
 g(x, y)y\leq -d y^2+exy
 \end{equation*}
 for some positive constants~$d$ and~$e$.
 \item[H$_3$] $L_g<\lambda_1$\,.
\item[H$_4$] $W$ is a Q-Wiener processes which has the following series expansion
\begin{equation}
W(t,x)=\sum_{i=1}^\infty \sqrt{q_i}e_i(x)\beta_i(t)
\end{equation}
and
\begin{equation*}
Qe_i=q_ie_i\,,\quad i=1,2,\ldots\,.
\end{equation*}
 Moreover, $ \operatorname{tr}Q<\infty$\,.
 \end{description}

By the above assumptions we have the following results on the averaging approximation to the slow-fast \spde{}~(\ref{e:ue})--(\ref{e:ve})~\cite{WR08}.
\begin{theorem} \label{thm:mixing}
Assume \textbf{H$_2$}.
For any fixed slow part $u\in H$\,, the fast system~(\ref{e:ve}) has a unique stationary solution,~$\eta^{\e,u}(t)$, with distribution~$\mu^u$ independent of~$\e$.
Moreover, the stationary measure~$\mu^u$ is exponentially mixing.
Also, $\eta^{\e, u}$~is differentiable with respect to~$u$ with Fr\'echet derivative
\begin{equation*}
D_u\eta^{\e,u}\leq D_v
\end{equation*}
for some positive constant~$D_v$ which is independent of~$\e$ and~$u$\,.
 \end{theorem}
Given the stationary measure~$\mu^u$ for the fast part, we define the following deterministic averaged equation
\begin{eqnarray}
 du&=&\big[A u+\bar{f}(u)\big]\,dt \,, \label{e:au1} \\
  u(0)&=&u_0 \,, \label{e:au2}
 \end{eqnarray}
where the average
\begin{equation}\label{e: fbar}
 \bar{f}(u)=\int_Hf(u, v)\mu^u(dv) \,.
 \end{equation}
Denote by~$\rho_{0T}$ the metric on space~$C(0, T; H)$  with
\begin{equation*}
\rho_{0T}(u,v)=\max_{0\leq t\leq T}\|u(t)-v(t)\|_0\quad \text{for all } u\,,v\in C(0, T; H)\,,
\end{equation*}
then we have the following theorem.
\begin{theorem} \label{thm:averaging}
 Assume \textbf{H$_1$}--\textbf{H$_4$}.
Given some $T>0$\,, for any $u_0\in H$\,, solutions~$u^\e(t, u_0)$ of~(\ref{e:ue}) converges in probability to~$u$ in~$C(0, T; H)$ which solves~(\ref{e:au1})--(\ref{e:au2}).
Moreover, the convergence rate is~$1/2$; that is, for any $\kappa>0$
\begin{equation*}
 \mathbb{P}\left\{\rho_{0T}(u^\e, u)\leq
 C^\kappa_T\sqrt{\e}\right\}>1-\kappa
\end{equation*}
for some positive constant $C^\kappa_T>0$\,.
 \end{theorem}

Now define the deviation between solutions~$u^\e$ and averaged solution~$u$,
\begin{equation}\label{e:ze}
 z^\e(t)=\frac{1}{\sqrt{\e}}(u^\e-u) \,.
\end{equation}
For $\e=1$ we write $\eta^{1,u}=\eta^u$\,.
Then we  have
 \begin{theorem}\label{thm:deviation}
 The deviation~$z^\e$ converges in distribution to a stochastic process~$z$ in the space~$C(0, T; H)$ which solves the \spde
 \begin{equation}\label{e:z}
  \dot{z}=Az+\overline{f'_u}(u)z+\sqrt{B(u)}\dot{\overline{W}}
 \end{equation}
 where $B(u): H\rightarrow H$ is Hilbert--Schmidt with
 \begin{align*}&
 B(u)=2\int_0^\infty\mathbb{E}\left[(f(u, \eta^u(t))-\bar{f}(u))\otimes(f(u,
 \eta^u(0))-\bar{f}(u))\right]dt
 \\&
\overline{f'_u}(u)=\int_Hf'_u(u, v)\mu^u(dv)
\end{align*}
and $\overline{W}(t)$ is an $H$-valued cylindrical Wiener process with covariance operator~$\operatorname{Id}_H$.
 \end{theorem}
\begin{remark}
Assume $\{e_i\}_{i=1}^\infty$ be a standard eigenbasis of~$H$, then
$B(u)$ has the following series form
\begin{align*}&
 B(u)=2\sum_{i,j=1}^\infty\int_0^\infty\mathbb{E}\Big[\left\langle f(u, \eta^u(t))-\bar{f}(u), e_i\right\rangle
 \left\langle f(u, \eta^u(0))-\bar{f}(u), e_j\right\rangle \Big]dte_i\otimes ej
\end{align*}
where $e_i\otimes e_j$ is the tensor product of $e_i$ and $e_j$\,.
\end{remark}

Then formally we write the following averaged equation plus
deviation up to errors of~$\mathcal{O}(\e)$ as
\begin{eqnarray}
 d\tilde{u}^\e&=&\big[A\tilde{u}^\e+\bar{f}(\tilde{u}^\e)\big]\,dt+\sqrt{\e}\sqrt{B(\tilde{u}^\e)}\,d\overline{W}(t) \,, \label{e:au1+deviation} \\
  \tilde{u}^\e(0)&=&u_0 \,. \label{e:au2+deviation}
 \end{eqnarray}

A family of random process~$\{u^\e\}_\e$ in space~$C(0, T; H)$
is said to satisfy the \ldp{} with rate function~$I$ if~\cite{FW98}
\begin{enumerate}
  \item(lower bound) for any $\phi\in C(0,T; H)$ and any $\delta,\gamma>0$\,, there is an~$\e_0>0$ such that
for~$0<\e<\e_0$
\begin{equation*}
 \mathbb{P}\left\{\rho_{0T}(u^\e, \phi)\leq \delta\right\}\geq  \exp\left\{-\frac{I(\phi)+\gamma}{\e}
\right\}\,,
\end{equation*}
\item(upper bound) for any $r>0$ and any $\delta,\gamma>0$\,,
there is a $\e_0>0$ such that for any $0< \e<\e_0$
\begin{equation*}
\mathbb{P}\left\{\rho_{0T}(u^\e, K_T(r))\geq \delta\right\} \leq
\exp\left\{-\frac{r-\gamma}{\e}\right\}\,.
\end{equation*}

\end{enumerate}
Here the level set $K_T(r):=\{\phi\in C(0,T; H): I(\phi)\leq r\}$\,.
If $K_T(r)$~is compact, then rate function~$I$ is called a good one.

There are many known results on the \ldp{} for those \spde{}s in the form of~(\ref{e:au1+deviation})--(\ref{e:au2+deviation}) as $\e\rightarrow 0$\,, as mentioned in  Section~\ref{sec:intro}.
Thus one may expect to derive \ldp{} for~$\{u^\e\}_\e$ from~(\ref{e:au1+deviation})--(\ref{e:au2+deviation}).
However, it is difficult to obtain an exponential approximation in probability between~$u^\e$ and~$\tilde{u}^\e$ as $\e\rightarrow 0$\,, which is needed to pass the \ldp{} of~$\{\tilde{u}^\e\}_\e$ to~$\{u^\e\}_\e$.
But here our result shows that the rate function for~$\{u^\e\}_\e$ as $\e\rightarrow 0$ is exactly that of~$\{\tilde{u}^\e\}_\e$.
This also shows system~(\ref{e:au1+deviation})--(\ref{e:au2+deviation}) is indeed an effective approximate  description of the macroscopic behaviour of the slow-fast system~(\ref{e:ue})--(\ref{e:ve}), even when one considers the exit problem caused by small noise perturbation for~$u^\e$ in the system of slow-fast \spde{}~(\ref{e:ue})--(\ref{e:ve}).

For our purposes we study the \ldp{} for a series of auxiliary slow-fast stochastic systems by a contraction principle and an approximation argument.
Then a controlled approximation derives the \ldp{} for~$\{u^\e\}_\e$.
We first give some contraction principles~\cite{DZ98} which are used in our approach.
\begin{lemma}\label{lem:contr-1}
Let $\mathcal{X}$ and~$\mathcal{Y}$ be Hausdorff topological spaces and $\Phi: \mathcal{X}\rightarrow \mathcal{Y}$ a continuous function.
Consider a good rate function $I: \mathcal{X}\rightarrow [0, \infty]$\,.
Then
\begin{enumerate}
  \item for each $y\in\mathcal{Y}$\,,
  \begin{equation*}
\tilde{I}(y)=\inf_{x\in\mathcal{X}}\{I(x): y=\Phi(x)\}
  \end{equation*}
is a good rate function on~$\mathcal{Y}$, where the infimum over an empty set is taken as~$\infty$;
  \item if $I$ controls the \ldp{} associated with a family of probability measures~$\{\mu_\e\}_\e$ on~$\mathcal{X}$, then $\tilde{I}$~controls the \ldp{} associated with the family of probability measures $\{\mu_\e \circ \Phi^{-1}\}_\e$ on~$\mathcal{Y}$.
\end{enumerate}
\end{lemma}

To introduce a generalized contraction principle we give the
following definitions.
\begin{definition}
Let $(\mathcal{Y}, d)$ be a metric space.
Then the probability measures $\{\mu_\e\}$ and~$\{\tilde{\mu}_\e\}_\e$ on~$\mathcal{Y}$ are called exponentially equivalent if there exists a probability space $\{\Omega, \mathcal{B}_\e, \mathbb{P}_\e\}$ and two families of $\mathcal{Y}$-valued random variables $\{Z_\e\}_\e$ and~$\{\tilde{Z}_\e\}_\e$ with joint laws $\{\mathbb{P}_\e\}_\e$ and marginals $\{\mu_\e\}$ and~$\{\tilde{\mu}_\e\}_\e$ respectively, such that for each $\delta>0$\,, the set $\{\omega: (\tilde{Z}_\e, Z_\e)\in \Gamma_\delta \}$ is $\mathcal{B}_\e$ measurable, and
\begin{equation*}
\limsup_{\e\rightarrow 0}\e\log\mathbb{P}_\e(\Gamma_\delta)=-\infty
\end{equation*}
where $\Gamma_\delta=\{(\tilde{y}, y):
d(\tilde{y},y)>\delta\}\subset \mathcal{Y}\times \mathcal{Y}$\,.
\end{definition}
As far as the \ldp{} is concerned, exponentially equivalent measures
are indistinguishable.
\begin{lemma}\label{lem:contr-2}
If an \ldp{} with good rate function $I(\cdot)$ holds for the probability measures $\{\mu_\e\}$, which are exponentially equivalent to $\{\tilde{\mu}_\e\}$, then the same \ldp{} holds for $\{\tilde{\mu}_\e\}$.
\end{lemma}

%Moreover we need the following generalized  contraction principle.
%\begin{lemma}\label{lem:contr-3}
%Let $(\mathcal{X}, d_1)$ and $(\mathcal{Y}, d_2)$ be Polish spaces
%and a family of probability measures $\{\mu_\e\}_\e$ satisfying
%\ldp{} with good rate function $I$\,. Suppose $\Phi_N:
%\mathcal{X}\rightarrow \mathcal{Y}$ is continuous mapping for any
%$N\in\mathbb{N}$ such that
%\begin{equation*}
%\lim_{N\rightarrow\infty}\limsup_{\e\rightarrow
%0}\log\mu_\e(d_2(\Phi_N, \Phi)>\delta)=-\infty\,,\quad  \text{for
%all} \quad \delta>0
%\end{equation*}
%where $\Phi:\mathcal{X}\rightarrow\mathcal{Y}$ is a measurable
%mapping, then there is a continuous function $\tilde{\Phi}:
%\{I<\infty\}\rightarrow \mathcal{Y}$ such that
%\begin{equation*}
%\lim_{N\rightarrow\infty}\sup_{I\leq r}d_2(\Phi_N,
%\tilde{\Phi})=0\,\quad \text{for all}\quad r>0
%\end{equation*}
%and $\{\mu_\e\circ \Phi^{-1}\}_\e$ satisfies \ldp{} with  good rate
%function $\tilde{I}$ with
%\begin{equation*}
%\tilde{I}(\phi):=\inf_{z\in\tilde{\Phi}^{-1}(\phi)}I(z)\,.
%\end{equation*}
%\end{lemma}
%
%

Furthermore, we need the following assumption
\begin{description}
  \item[H$_5$] There is a positive constant $c_0$ and $c_1$ such that
\begin{equation*}
 \langle B(\phi)h, h\rangle \geq c_0\|h\|_0^2\quad \text{and}
 \quad \langle D B(\phi) h, h\rangle<c_1\|h\|_0^2\,,\quad \text{for all } \phi\,, h\in H
\end{equation*}
\end{description}
where $D B$ is the Fr\'echlet derivative of $B$\,.

\begin{remark}
Under the above assumption, $\sqrt{B(\phi)}$ is Lipschitz continuous
in~$\phi$. The following is one simple example of $f$ such
that $B(\phi)$ satisfies assumption $\textbf{H}_5$\,,
\begin{equation*}
f(u,v)=f_1(u)+f_2(v)
\end{equation*}
with $f_1$ and $f_2$ both Lipschitz continuous. If the
stationary Gaussian process~$\eta^u=\eta$ is independent of~$u$,  by
the expression of~$B(u)$ in Theorem~\ref{thm:deviation},
\begin{align*}&
 B(u)=2\int_0^\infty\mathbb{E}\Big[(f_2(\eta(t))-\mathbb{E}f_2(\eta))
 \otimes(f_2(\eta(0))-\mathbb{E}f_2(\eta))\Big]dt\,.
\end{align*}
Then  $B(u)$, which is independent of~$u$, satisfies
assumption~$\textbf{H}_5$.
Section~\ref{sec:esfsrd} details one example where $f_2(v)$ is linear in the fast~$v$.
\end{remark}

Now for the slow-fast stochastic system~(\ref{e:ue}) we define the
following skeleton equation
\begin{equation}\label{e:skeleton}
\dot{\phi}=A\phi+\bar{f}(\phi)+\sqrt{B(\phi)}h\,, \quad
\phi(0)=u_0\,.
\end{equation}
And define the rate function
\begin{equation}\label{e:rate function}
I_u(\phi)=\inf_{h\in L^2(0, T;
H)}\left\{\frac12\int_0^T\|h(s)\|_0^2\,ds: \phi=\phi^{ h} \right\}
\end{equation}
where $\phi^{h}$ solves~(\ref{e:skeleton}) and with
$\inf{\emptyset}=+\infty$\,. Then we give our main result.
\begin{theorem}\label{thm:main}
Assume $\textbf{H}_1$--$\textbf{H}_5$\,. For any $T>0$\,,
\{$u^\e\}_\e$ satisfies the \ldp{} with good rate function~$I_u$ on
space~$C(0,T; H)$.
\end{theorem}

\begin{remark}
By the large deviations principle for stochastic evolutionary
equations~\cite{PZ92}, under assumptions
$\textbf{H}_1$--$\textbf{H}_5$\,, $\{\tilde{u}^\e\}_\e$ satisfies the \ldp{} with rate function $I_u(\phi)$ in space $C(0, T; H)$ for any $T>0$\,. So the averaged equation plus
deviation~(\ref{e:au1+deviation})--(\ref{e:au2+deviation}) predicts the metastability of~(\ref{e:ue})--(\ref{e:ve}).
\end{remark}

\section{Exponential tightness}\label{sec:exp tight}
We prove some exponential tightness results which are crucial for
the \ldp{} estimates.

For any $\psi\in C(0, T; H)$ consider the following system
\begin{eqnarray}\label{e:tilde-u-psi}
d \tilde{u}^{\e, \psi}&=&\left[A\tilde{u}^{\e, \psi}+f(\psi, \tilde{v}^{\e, \psi})\right]\,dt\,,\quad \tilde{u}^{\e, \psi}(0)=u_0\,,\\
d\tilde{v}^{\e,\psi}&=&\frac{1}{\e}\left[A\tilde{v}^{\e,\psi}+g(\psi,
\tilde{v}^{\e,\psi})\right]\,dt+\frac{1}{\sqrt{\e}}dW(t)\,,\quad
\tilde{v}^{\e,\psi}(0)=v_0\,.\label{e:tilde-v-psi}
\end{eqnarray}
Then we have the following exponential tightness result.
\begin{lemma}\label{lem:exp-tight-1}
Fix $\psi\in C(0, T; H)$.  For any $T>0$ and $\e_0>0$\,,
$\{\tilde{u}^{\e, \psi}\}_{0\leq \e\leq \e_0}$ is exponentially
tight in space $C(0, T; H)$.
\end{lemma}
\begin{proof}
This result follows via a uniform estimate to $\tilde{u}^{\e,
\psi}$ in space $C(0, T; H_0^1)\cap C^\alpha(0, T; H)$ for some
$\alpha>0$\,. By~(\ref{e:tilde-u-psi})
\begin{equation*}
\tilde{u}^{\e,
\psi}(t)=u_0+\int_0^tA\tilde{u}^{\e,\psi}(s)\,ds+\int_0^tf(\psi(s),
\tilde{v}^{\e, \psi}(s))\,ds\,.
\end{equation*}
Then by the increasing  property of $f$
\begin{eqnarray*}
&&\|\tilde{u}^{\e,\psi}(t)-\tilde{u}^{\e, \psi}(\tau)\|^2_0\\&\leq&
2\left[\int_\tau^t\|A\tilde{u}^{\e,
\psi}(s)\|_0\,ds\right]^2+2\left[\int_\tau^t\|f(\psi(s),
\tilde{v}^{\e, \psi}(s))\|_0\,ds\right]^2\\
&\leq&2|t-\tau|\int_0^T\|A\tilde{u}^{\e,\psi}(s)\|_0^2\,ds+
2|t-\tau|\int_\tau^t\|f(\psi(s), \tilde{v}^{\e, \psi}(s))\|_0^2\,ds
\\&\leq&2|t-\tau|\int_0^T\|A\tilde{u}^{\e,\psi}(s)\|_0^2\,ds+
2|t-\tau|\int_0^T\left[a\|\psi(s)\|_0^2+\|\tilde{v}^{\e,
\psi}(s)\|_0^2+c\right]ds\,,
\end{eqnarray*}
and
\begin{eqnarray*}
\frac12\frac{d}{dt}\|\tilde{u}^{\e,
\psi}\|_1^2&=&-\|A\tilde{u}^{\e,\psi}\|_0^2-\langle f(\psi,
\tilde{v}^{\e, \psi}), A\tilde{u}^{\e, \psi}\rangle\\
&\leq&-\frac12\|A\tilde{u}^{\e,\psi}\|_0^2+\frac12\|f(\psi,
\tilde{v}^{\e, \psi})\|_0^2\\
&\leq&-\frac12\|A\tilde{u}^{\e,\psi}\|_0^2+\frac12\left[a\|\psi\|_0^2+\|\tilde{v}^{\e,
\psi}\|_0^2+c\right].
\end{eqnarray*}
That is,
\begin{eqnarray*}
&&\int_0^T\|A\tilde{u}^{\e, \psi}(s)\|_0^2\,ds+\sup_{0\leq s\leq
T}\|\tilde{u}^{\e, \psi}(s)\|_1^2\\&\leq&
\|u_0\|_1^2+a\int_0^T\|\psi(s)\|_0^2\,ds+\int_0^T\|\tilde{v}^{\e,
\psi}(s)\|_0^2\,ds+cT\,.
\end{eqnarray*}
Then to give a uniform estimate in space $C(0,T; H_0^1)\cap
C^\alpha(0,T; H)$ for some $\alpha>0$\,, we just need to estimate
$\int_0^T\|\tilde{v}^{\e, \psi}(s)\|_0^2\,ds$\,.

Next we estimate $\int_0^T\|\tilde{v}^{\e,
\psi}\|_1^2\,ds$\,. The Lipschitz property of $g$ and applying the
It\^o formula to $\|\tilde{v}^{\e, \psi}\|^2_0$ yield for some
positive constants $c'$ and~$c''$
\begin{eqnarray*}
\frac12\frac{d}{dt}\|\tilde{v}^{\e,
\psi}\|_0^2&=&-\frac{1}{\e}\|\tilde{v}^{\e,
\psi}\|_1^2+\frac{1}{\e}\langle g(\psi, \tilde{v}^{\e, \psi}),
\tilde{v}^{\e, \psi}\rangle+\frac{1}{\sqrt{\e}}\langle
\tilde{v}^{\e,\psi},\dot{W}\rangle +\frac{1}{2\e}\operatorname{tr} Q\\
&\leq& -\frac{c'}{\e}\|\tilde{v}^{\e,
\psi}\|_1^2+\frac{c''}{\e}\|\psi\|_0^2+\frac{1}{2\e}\operatorname{tr}Q+
\frac{1}{\sqrt{\e}}\langle \tilde{v}^{\e,\phi}, \dot{W}\rangle\,.
\end{eqnarray*}
Then
\begin{eqnarray*}
&&\e\|\tilde{v}^{\e, \psi}(t)\|_0^2+2c'\int_0^t\|\tilde{v}^{\e,
\phi}(s)\|_1^2\,ds-\e\|v_0\|_0^2\\
&\leq &2c''\int_0^t\|\psi(s)\|_0^2\,ds+\operatorname{tr}Q
t+2\sqrt{\e}\int_0^t\langle \tilde{v}^{\e, \psi}(s), dW(s)
\rangle\,.
\end{eqnarray*}
Now define $M^\e_t=\int_0^t\langle \tilde{v}^{\e, \psi}(s), dW(s)
\rangle$ and $\lambda_0=c'\lambda_1/2q_{\max}$ with $q_{\max}=\max_i
q_i$\,. And denote by $\langle M\rangle^\e_t$ the covariance of
$M_t^\e$. Then we have
\begin{eqnarray*}
c'\int_0^t\|\tilde{v}^{\e, \psi}(s)\|_1^2\,ds&\leq&
\|v_0\|_0^2+2c''\int_0^t\|\psi(s)\|_0^2\,ds+\operatorname{tr}Q
t+\sqrt{\e}M_t^\e-\frac{\lambda_0}{2}\langle M \rangle_t^\e
\\&&{}+
\frac{\lambda_0}{2}\langle M \rangle_t^\e-c'\int_0^t\|\tilde{v}^{\e,
\psi}(s)\|_1^2\,ds\\
&\leq&\|v_0\|_0^2+2c''\int_0^t\|\psi(s)\|_0^2\,ds+\operatorname{tr}Q
t+\sqrt{\e}M_t^\e-\frac{\lambda_0}{2}\langle M \rangle_t^\e\,.
\end{eqnarray*}
And by the exponential martingale inequality we have
\begin{equation*}
\mathbb{P}\left\{\sqrt{\e}M_t^\e-\frac{\lambda_0}{2}\langle M
\rangle_t^\e>\delta\right\}
=\mathbb{P}\left\{\frac{\lambda_0}{\sqrt{\e}}M_t^\e-\frac{\lambda^2_0}{2\e}\langle
M \rangle_t^\e>\frac{\delta\lambda_0}{\e}\right\}
\leq e^{-{\delta\lambda_0}/{\e}}
\end{equation*}
which yields the exponential tightness of
$\{\tilde{u}^{\e,\psi}\}_\e$\,.
\end{proof}
Similar method yields the following result.
\begin{lemma}\label{lem:exp-tight-2}
For any $T>0$\,, $\{u^\e\}_\e$ is exponentially tight in space $C(0,
T; H)$.
\end{lemma}
\begin{proof}
Follow the proof of Lemma~\ref{lem:exp-tight-1}, we just need to
 estimate
$\int_0^T\|u^\e(s)\|_0^2\,ds+\int_0^T\|v^\e(s)\|_0^2\,ds$\,. By our assumption we  estimate
\begin{equation*}
\int_0^T\|u^\e(s)\|_1^2\,ds+\int_0^T\|v^\e(s)\|_1^2\,ds\,.
\end{equation*}
For this, applying the It\^o formula to~$\|v^\e\|_0^2$ and by the
assumption on~$g$ we have
\begin{equation}\label{e:1}
\frac{1}{2}\frac{d}{dt}\|v^\e\|_0^2\leq
-\frac{1}{\e}\|v^\e\|_1^2-\frac{d}{\e}\|v^\e\|_0^2+\frac{e}{\e}\langle
u^\e, v^\e\rangle+\frac{1}{2\e}\operatorname{tr}Q+\frac{1}{\sqrt{\e}}\langle
v^\e, \dot{W} \rangle\,.
\end{equation}
Moreover,
\begin{equation}\label{e:2}
\frac12\frac{d}{dt}\|u^\e\|_0^2\leq
-\|u^\e\|_1^2-a\|u^\e\|_0^2+b\langle u^\e, v^\e \rangle+c \,.
\end{equation}
Then $\e(\ref{e:1})+(\ref{e:2})$ yields
\begin{eqnarray*}
&&2\int_0^T\|u^\e(s)\|_1^2\,ds+2\int_0^T\|v^\e(s)\|_1^2\,ds
\\&\leq& \|u_0\|_0^2+\|v_0\|_0^2+C_T+\operatorname{tr}QT+
\sqrt{\e}\int_0^t\langle v^\e(s), dW(s)\rangle\,.
\end{eqnarray*}
Then a similar argument as for Lemma~\ref{lem:exp-tight-1} yields the
exponential tightness. The proof is complete.
\end{proof}

\section{LDP for some auxiliary systems}
\label{sec:aux-sys}
Next we study the \ldp{} for some auxiliary systems.

First, let $\psi\in H$\,, which is independent of time, and consider the
following system
\begin{eqnarray}\label{e:u-psi}
d u^{\e, \psi}&=&\left[Au^{\e, \psi}+f(\psi, v^{\e, \psi})\right]\,dt\,,\quad u^{\e, \psi}(0)=u_0\,,\\
dv^{\e,\psi}&=&\frac{1}{\e}\left[Av^{\e,\psi}+g(\psi,
v^{\e,\psi})\right]\,dt+\frac{1}{\sqrt{\e}}dW(t)\,,\quad
v^{\e,\psi}(0)=v_0\,.\label{e:v-psi}
\end{eqnarray}
We study the \ldp{} for~$\{u^{\e,\psi}\}_\e$. For this define
\begin{equation*}
\xi^{\e,\psi}(t)=\int_0^t\left[f(\psi,
v^{\e,\psi}(s))-\bar{f}(\psi)\right]\,ds\,.
\end{equation*}
Then
\begin{equation}\label{e:u-e-psi}
\dot{u}^{\e,\psi}=Au^{\e, \psi}+\bar{f}(\psi)+\dot{\xi}^{\e,
\psi}\,,\quad u^{\e,\psi}(0)=u_0\,.
\end{equation}
Now in order to obtain the \ldp{} for~$\{u^{\e,\psi}\}_\e$ we study
the \ldp{} for~$\{\bar{u}^{\e, \psi}\}_\e$ which solves
\begin{equation}\label{e:bar-u}
\dot{\bar{u}}^{\e, \psi}=A\bar{u}^{\e,
\psi}+\bar{f}(\psi)+\dot{\bar{\xi}}^{\e,\psi}\,,\quad \bar{u}^{\e,
\psi}(0)=u_0
\end{equation}
with
% follow an
%approximation discussion~\cite{RWW}. For any fixed $N\geq 1$\,, let
%$H_N:=\text{span}\{e_i: i=1, \ldots, N\}$~and~$P_N: H\rightarrow
%H_N$ be the orthogonal projection. Let $u^{\e, \psi}_N$ be the
%solution of
%\begin{equation}\label{e:ueN}
%\dot{u}^{\e, \psi}_N=Au^{\e, \psi}_N+\bar{f}(u^{\e,
%\psi}_N)+P_N\dot{\bar{\xi}}^{\e, \psi}\,,\quad u^{\e, \psi}_N(0)=u_0
%\end{equation}
%with
\begin{eqnarray*}
\bar{\xi}^{\e,\psi}(t)=\int_0^t\left[f(\psi,
\eta^{\e,\psi}(s))-\bar{f}(\psi)\right]\,ds
\end{eqnarray*}
where $\eta^{\e,\psi}$ is defined in Theorem~\ref{thm:mixing}\,. Now
define the  skeleton equation to~(\ref{e:bar-u}) with~$\psi\in H$
\begin{equation}\label{e:skeleton1}
\dot{\phi}^{\psi, h}=A\phi^{\psi,h}+\bar{f}(\psi)+\sqrt{B(\psi)}h\,,
\quad \phi^{\psi, h}(0)=u_0
\end{equation}
for $h\in L^2(0, T; H)$\,. Equation (\ref{e:skeleton1}) is a linear
 equation, so for any~$\psi\in H$ and~$h\in
L^2(0,T; H)$  there is a unique~$\phi^{\psi, h}\in C(0,T; H)$
solving~(\ref{e:skeleton1}). Then define the functional
\begin{equation*}
I_u^\psi(\phi)=\inf_{h\in L^2(0, T;
H)}\left\{\frac12\int_0^T\|h(s)\|_0^2\,ds: \phi=\phi^{\psi,
h}\right\}.
\end{equation*}
We prove the following \ldp{} result for
$\{\bar{u}^{\e,\psi}\}_\e$\,.
\begin{theorem}\label{thm:ldp-bar-u}
Fix $\psi\in H$\,.  For any $T>0$\,, $\{\bar{u}^{\e, \psi}\}_\e$
satisfies \ldp{} in space $C(0, T; H)$ with a good rate function
$I_u^\psi$. Moreover~$I_u^\psi$ is lower semicontinuous in
$\psi\in H$\,.
\end{theorem}

\begin{proof}
First by the same discussion in the proof to Lemma
\ref{lem:exp-tight-1}\,, for any~$T>0$\,, $\{\bar{u}^{\e,
\psi}\}_\e$ is exponential tight in space~$C(0, T; H)$. Then there
exists~$\{K_R\}_R$ which is a nondecreasing family of compact sets
such that
\begin{equation}\label{e:K-R}
\mathbb{P}\{\bar{u}^{\e,\psi}\in K_R\} \geq 1-e^{-{R}/{\e}}\,.
\end{equation}
Now for any fixed $N\geq 1$\,, let $H_N:=\text{span}\{e_i: i=1,
\ldots, N\}$ and~$P_N: H\rightarrow H_N$ be the orthogonal
projection. Then~$\bar{\xi}_{N}^{\e, \psi}:=P_N\bar{\xi}^{\e,
\psi}\in C(0, T; H_N)$ satisfies \ldp{} with good rate
function~\cite{FW98}
\begin{equation*}
I_{\xi,
N}^\psi(\phi_N)=\inf_{h_N}\left\{\frac12\int_0^T\|h_N(t)\|_{H_N}^2\,dt:
\sqrt{B_N(\psi)}h_N=\phi_N\right\}
\end{equation*}
where $\phi_N\in C(0,T; H_N)$ and
\begin{equation*}
B_N(\psi)=2\int_0^\infty\mathbb{E}\left[(P_Nf(\psi,\eta^\psi(t))-P_N\bar{f}(\psi))\otimes(P_Nf(\psi,
\eta^\psi(0))-P_N\bar{f}(\psi)) \right]\,dt\,.
\end{equation*}
Introduce process $\bar{u}_N^{\e,\psi}$  solving
\begin{equation*}
\dot{\bar{u}}^{\e, \psi}_N=A\bar{u}^{\e,
\psi}_N+\bar{f}_N(\psi)+\dot{\bar{\xi}}^{\e, \psi}_N\,, \quad
\bar{u}^{\e,\psi}_N(0)=P_Nu_0
\end{equation*}
with $\bar{f}_N=P_N\bar{f}$\,.
 By the continuity of the map $\bar{\xi}^{\e, \psi}_N\mapsto
\bar{u}^{\e, \psi}_N$ in space $C(0,T; H)$ and the contraction
principle Lemma~\ref{lem:contr-1}\,, $\{\bar{u}^{\e, \psi}_N\}_\e$
satisfies \ldp{} with a good rate function
\begin{equation*}
I_{u, N}^\psi(\phi)=\inf_{h\in L^2(0,T;
H)}\left\{\frac12\int_0^T\|h(s)\|_0^2ds: \phi_N=\phi_N^{\psi,h}
\right\}
\end{equation*}
where $\phi^{\psi, h}\in C(0,T; H)$ solves the following equation
\begin{equation*}
\dot{\phi}^{\psi,h}_N=A\phi_N^{ \psi,
h}+\bar{f}_N(\psi)+\sqrt{B_N(\psi)}h\,,\quad \phi^{\psi,
h}(0)=P_Nu_0\,.
\end{equation*}
Moreover $\phi_N^{\psi, h}\rightarrow \phi^{\psi, h}$ as
$N\rightarrow\infty$ in space $C(0, T;H)$\,. Then we have
\begin{equation}
I_u^\psi(\phi)=\lim_{N\rightarrow\infty}I_{u,N}^\psi(P_N\phi)=\sup_NI_{u,N}^\psi(P_N\phi)\,,
\quad \phi\in C(0, T; H)
\end{equation}
which is a good rate function.

Now for any $\gamma$\,, $\delta>0$\,, $\phi\in C(0,T; H)$ and
$R>0$\,, for $\bar{u}^{\e,\psi}\in K_R$\,, there is~$N_R(\delta,
\phi)>0$ such that
\begin{equation*}
\rho_{0T}(\bar{u}^{\e, \psi}, \bar{u}_N^{\e, \psi})+\rho_{0T}(\phi,
P_N\phi)<\delta/2\,, \quad N>N_R(\delta, \phi)\,.
\end{equation*}
Then there is $\e_0>0$ such that for $0<\e<\e_0$
\begin{eqnarray}
&&\mathbb{P}\left\{\rho_{0T}(\bar{u}^{\e, \psi}, \phi)<\delta\right\}\\
&\geq& \mathbb{P}\left\{\rho_{0T}(\bar{u}^{\e, \psi}_N,
P_N\phi)<\delta/2\,,\quad N>N_R(\delta, \phi)  \;\Big|\;  \bar{u}^{\e,
\psi}\in
K_R\right\}\nonumber\\
&\geq& \mathbb{P}\left\{\rho_{0T}(\bar{u}^{\e, \psi}_N,
P_N\phi)<\delta/2  \right\}\mathbb{P}\left\{\bar{u}^{\e, \psi}\in
K_R\;\Big|\;  \rho_{0T}(\bar{u}^{\e, \psi}_N, P_N\phi)<\delta/2
\right\}\nonumber\\
&\geq&\exp\left\{-\frac{I_{u,N}^\psi(P_N\phi)+\gamma}{\e}\right\}\mathbb{P}\left\{\bar{u}^{\e,
\psi}\in K_R\;\Big|\;  \rho_{0T}(\bar{u}^{\e, \psi}_N, P_N\phi)<\delta/2
\right\}\nonumber\\
&\geq&\exp\left\{-\frac{I_{u}^\psi(\phi)+\gamma}{\e}\right\}\mathbb{P}\left\{\bar{u}^{\e,
\psi}\in K_R\;\Big|\;  \rho_{0T}(\bar{u}^{\e, \psi}_N, P_N\phi)<\delta/2
\right\}\,.\nonumber
\end{eqnarray}
Passing to the limit $R\rightarrow\infty$ and noticing that
$\lim_{R\rightarrow\infty}\mathbb{P}\{\bar{u}^{\e,\psi}\in
K_R\}=1$\,, the above estimate yields the lower bound estimate.

Now for any $r>0$ and $\gamma>0$\,, there is an $\e_1>0$  and $R>0$
such that for any $0<\e<\e_1$
\begin{equation}\label{e:K-c}
\mathbb{P}\left\{\bar{u}^{\e, \psi}\in K^c_R \right\}\leq
\exp\left\{-\frac{r-\gamma}{\e} \right\}\,.
\end{equation}
Let $K_T(r)=\{\phi\in C(0, T; H): I_u^\psi(\phi)\leq r\}$ and $K_{T,
N}(r)=\{\phi_N\in C(0, T;H_N): I_{u,N}^\psi(\phi_N)\leq r\}$\,. Then
for $\bar{u}^{\e, \psi}\in K_R$\,, there is $N(R, r)<0$ such that
for~$N>N(R, r)$
\begin{eqnarray}
&&\mathbb{P}\left\{\rho_{0T}(\bar{u}^{\e, \psi},
K_T(r))>\delta\;\big|\;\bar{u}^{\e, \psi}\in
K_R\right\}\nonumber\\&\leq&
\mathbb{P}\left\{\rho_{0T}(\bar{u}_N^{\e, \psi},
K_{T,N}(r))>\delta/2\;\big|\;\bar{u}^{\e, \psi}\in K_R\right\}\nonumber\\
&\leq& \exp\left\{-\frac{r-\gamma}{\e} \right\}\label{e:K}
\end{eqnarray}
for $\e$ is small enough.  Then (\ref{e:K-c})--(\ref{e:K}) yield the
upper bound estimate.

Next we prove the second result. This is followed by proving the set
\begin{equation*}
\Psi(r):=\left\{\psi\in H: I_u^\psi(\phi)\leq r \right\}
\end{equation*}
is closed for any $r>0$ and $\phi\in C(0, T;H)$.  For this proof, let
$\{\psi_n\}_n\subset\Psi(r)$~and~$\psi_n\rightarrow \psi$ in space~$H$. Moreover, there is a sequence $h_n\in L^2(0,T; H)$ such
that~$\phi^{\psi_n,h_n}=\phi$ with
\begin{equation*}
\dot{\phi}^{\psi_n,
h_n}=A\phi^{\psi_n,h_n}+\bar{f}(\psi_n)+\sqrt{B(\psi_n)}h_n\,,\quad
\phi^{\psi_n,h_n}(0)=u_0
\end{equation*}
and
\begin{equation*}
\frac12\int_0^T\|h_n(s)\|_0^2\,ds<r+\frac{1}{n}\,.
\end{equation*}
Then there is a subsequence, which we still denote by $h_n$\,,
weakly convergent to some $h\in L^2(0, T; H)$\,. We show that
$\phi^{\psi, h}=\phi$ with
\begin{equation*}
\dot{\phi}^{\psi,
h}=A\phi^{\psi,h}+\bar{f}(\psi)+\sqrt{B(\psi)}h\,,\quad
\phi^{\psi,h}(0)=u_0
\end{equation*}
and
\begin{equation*}
\frac12\int_0^T\|h(s)\|_0^2\,ds\leq r\,.
\end{equation*}
By the assumptions $\textbf{H}_1$ and $\textbf{H}_5$, we have
\begin{equation*}
\bar{f}(\psi_n)\rightarrow \bar{f}(\psi) \quad \text{in}\quad H
\end{equation*}
and
\begin{equation*}
\sqrt{B(\psi_n)}\rightarrow \sqrt{B(\psi)} \quad \text{in}\quad
\mathcal{L}(H)\,.
\end{equation*}
Then
\begin{equation*}
\sqrt{B(\psi_n)}h_n\rightarrow \sqrt{B(\psi)}h \quad \text{weakly
in}\quad L^2(0, T; H)\,.
\end{equation*}
The proof is complete.

\end{proof}

By the above result we have the following corollary.
\begin{corollary}\label{cor:ldp-u-psi}
For fixed $\psi\in H$\,,  $\{u^{\e, \psi}\}_\e$ satisfies \ldp{}
with good rate function $I_u^\psi$ in space $C(0, T; H)$ for any
$T>0$\,.
\end{corollary}
\begin{proof}
We show that $u^{\e, \psi}$ is exponentially equivalent in
probability to $\bar{u}^{\e,\psi}$\,. Let
$U^{\e,\psi}=u^{\e,\psi}-\bar{u}^{\e,\psi}$\,, then
\begin{equation*}
\dot{U}^{\e,\psi}=AU^{\e,\psi}+\dot{\xi}^{\e,\psi}-\dot{\bar{\xi}}^{\e,\psi}\,,\quad
U^{\e,\psi}(0)=0\,.
\end{equation*}
Notice that
\begin{equation*}
\xi^{\e, \psi}(t)-\bar{\xi}^{\e,\psi}(t)=\int_0^t\left[f(\psi,
v^{\e,\psi}(s))-f(\psi,\eta^{\e, \psi}(s))\right]\,ds\,.
\end{equation*}
Then for any $T>0$
\begin{eqnarray}
\sup_{0\leq t\leq T}\|u^{\e, \psi}(t)-\bar{u}^{\e,\psi}(t)\|^2_0
&\leq&C_{\lambda_1}\sup_{0\leq t\leq T}\int_0^t\|f(\psi,v^{\e,
\psi}(s))-f(\psi,
\eta^{\e, \psi}(s))\|_0\,ds\nonumber \\
&\leq & L_fC_{\lambda_1}\int_0^T\|v^{\e, \psi}(s)-\eta^{\e,
\psi}(s)\|_0\,ds \label{e:xi-xi_bar}
\end{eqnarray}
for some positive constant $C_{\lambda_1}$ which depends on
$\lambda_1$\,. Let $\zeta^{\e, \psi}=v^{\e, \psi}-\eta^{\e,
\psi}$\,,
\begin{equation*}
\dot{\zeta}^{\e, \psi}=\frac{1}{\e}A\zeta^{\e,
\psi}+\frac{1}{\e}\left[g(\psi, v^{\e, \psi})-g(\psi, \eta^{\e,
\psi})\right]\,, \quad \zeta^{\e, \psi}(0)=v_0-\eta^{\e, \psi}(0)\,.
\end{equation*}
Then
\begin{equation*}
\frac12\frac{d}{dt}\|\zeta^{\e, \psi}(t)\|_0^2\leq
-\frac{\lambda_1-L_g}{\e}\|\zeta^{\e, \psi}(t)\|_0^2\,.
\end{equation*}
By assumption $\mathbf{H}_3$ we have
\begin{equation*}
\|\zeta^{\e, \psi}(t)\|_0^2\leq e^{-2tc/\e}\|v_0-\eta^{\e,
\psi}(0)\|_0^2
\end{equation*}
for some positive constant $c$ which is independent $\e$\,. Then by
(\ref{e:xi-xi_bar}) we have for any $T>0$
\begin{eqnarray*}
\sup_{0\leq t \leq T}\|u^{\e, \psi}(t)-\bar{u}^{\e,\psi}(t)\|_0
&\leq& L_fC_{\lambda_1}\int_0^Te^{-sc/\e}\|v_0-\eta^{\e, \psi}(0)\|_0\,ds\\
&\leq & \e C\|v_0-\eta^{\e, \psi}\|_0
\end{eqnarray*}
for some positive constant $C$ which is independent of $\e$\,.

Then we have  for any $\delta>0$\,,
\begin{equation*}
\mathbb{P}\left\{\sup_{0\leq t\leq
T}\|u^{\e,\psi}(t)-\bar{u}^{\e,\psi}(t)\|_0> \delta\right\}\leq
\mathbb{P}\left\{\e C\|v_0-\eta^{\e,\psi}(0)\|>\delta\right\}\,.
\end{equation*}
And by the Guassian property of $\eta^{\e, \psi}(0)$\,, we have for
any $\delta>0$
\begin{equation*}
\limsup_{\e\rightarrow 0}\e\ln\mathbb{P}\left\{\sup_{0\leq t\leq
T}\|u^{\e,\psi}(t)-\bar{u}^{\e,
\psi}(t)\|_0>\delta\right\}=-\infty\,.
\end{equation*}
Then the generalized contraction principle Lemma \ref{lem:contr-2}
completes the proof.\\
\end{proof}

Now we consider the special case that $\psi$~is a step function
which we denote by~$\psi^n$. Let $(\tilde{u}^{\e, \psi^n},
\tilde{v}^{\e, \psi^n})$
solve~(\ref{e:tilde-u-psi})--(\ref{e:tilde-v-psi}) with~$\psi$
replaced by~$\psi^n$, we show that $\{\tilde{u}^{\e, \psi^n}\}_\e$
satisfies the \ldp{} with a good rate function. Typically we choose
\begin{equation}\label{e:step-function}
\psi^n(t)=\sum_{i=0}^{n-1}\psi^n_i\chi_{[t_i, t_{i+1}]}
\end{equation}
with $\psi_i^n\in H$ and  $t_i=iT/n$\,, $i=0, 1,\ldots, n-1$\,. Let
$\tilde{u}^{\e, \psi^n_i}$ and~$\tilde{v}^{\e, \psi_i^n}$ be the
restriction of~$\tilde{u}^{\e,\psi^n}$ and~$\tilde{v}^{\e, \psi^n}$
on the time interval~$[t_i, t_{i+1}]$, respectively, then we have
\begin{equation*}
\dot{\tilde{u}}^{\e, \psi_i^n}=A\tilde{u}^{\e, \psi_i^n}+f(\psi_i^n,
\tilde{v}^{\e,\psi_i^n})
\end{equation*}
which satisfies the \ldp{} in space~$C([t_i, t_{i+1}], H)$ with rate
function~$I_u^{\psi_i^n}$ by Corollary~\ref{cor:ldp-u-psi}. Then we
have the following result

\begin{theorem}\label{thm:LDP-2}
For any $n$ and step function $\psi^n\in C(0, T; H)$ in the
form~(\ref{e:step-function}), $\{\tilde{u}^{\e, \psi^n}\}$
satisfies \ldp{} with good rate function
\begin{equation*}
I_u^{\psi^n}(\phi)=\sum_{i=0}^{n-1}I_u^{\psi^n_i}=\inf_{h\in
L^2(0,T; H)}\left\{\frac12\int_0^T\|h(s)\|_0^2\,ds:
\phi=\phi^{\psi^n, h} \right\}
\end{equation*}
where $\phi^{\psi^n, h}$ solves (\ref{e:skeleton1}) with $\psi$
replaced by step function $\psi^n$\,.
\end{theorem}
\begin{proof}
For any $\phi\in C(0, T; H)$ with $I_u^{\psi^n}(\phi)<\infty$\,,
let~$\phi_i(t)=\phi(t)\chi_{[t_i,\, t_{i+1}]}\in C([t_i,t_{i+1}];
H])$\,, then for any~$\delta, \gamma>0$\,, there is an $\e_0>0$
such that for any $0<\e<\e_0$\,,
\begin{equation*}
\mathbb{P}\left\{\max_{t_i\leq t\leq t_{i+1}}\|\tilde{u}^{\e,
\psi_i^n}(t)-\phi_i(t)\|_0\leq \delta
\right\}\geq\exp\left[-\frac{I_u^{\psi_i^n}(\phi_i)+\frac{\gamma}{n}}{\e}
\right].
\end{equation*}
Moreover, by the assumption on~$f$, for $n$~large enough for any
$\delta>0$
\begin{equation*}
\left\{\max_{t_i\leq t\leq t_{i+1}}\|\tilde{u}^{\e,
\psi^n_i}(t)-\psi_i(t)\|_0\leq\delta\right\}\supset\left\{\|\tilde{u}^{\e,
\psi^n_i}(t_i)-\psi_i(t_i)\|_0\leq\delta'(\delta) \right\}:=A_i
\end{equation*}
for some $\delta'(\delta)$ small enough. Now for fixed time~$t_i$,
$\tilde{u}^{\e, \psi_i^n}(t_i)$ is an $H$-valued random variable.
Then for small enough $\delta'=\delta'(\delta)$
\begin{eqnarray*}
&&\mathbb{P}\left\{\rho_{0T}(\tilde{u}^{\e,
\psi^n}, \phi)\leq\delta\right\}\\
&=& \mathbb{P}\left\{\max_{1\leq i\leq n,\, t_i\leq t\leq
t_{i+1}}\|\tilde{u}^{\e, \psi_i^n}(t)-\phi_i(t)\|_0\leq \delta \right\}\\
&\geq& \mathbb{P}\left\{\max_{1\leq i\leq n}\|\tilde{u}^{\e,
\psi_i^n}(t_i)-\phi_i(t_i)\|_0\leq \delta' \right\}
=\mathbb{P}\left\{A_1A_2\ldots A_n\right\}\\
&=&\mathbb{P}\left\{A_1\right\} \mathbb{P}\left\{A_2\mid A_1\right\}
\mathbb{P}\left\{A_3\mid A_1A_2\right\}
\cdots\mathbb{P}\left\{A_n\mid A_1A_2\ldots A_{n-1}\right\}\\
&\geq&\mathbb{P}\left\{ \max_{t_0\leq t\leq t_1}\|\tilde{u}^{\e,
\psi_0^n}(t)-\phi_0(t)\|_0\leq \delta'\right\}
\\&&{}\times
\mathbb{P}\left\{\max_{t_1\leq t\leq
t_2}\|\tilde{u}^{\e, \psi_1^n}(t)-\phi_1(t)\|_0\leq \delta'
\mid A_1\right\}
\\&&{}
\times\cdots
\times\mathbb{P}\left\{\max_{t_{n-1}\leq t\leq
t_n}\|\tilde{u}^{\e, \psi_{n-1}^n}(t)-\phi_{n-1}(t)\|_0\leq \delta'
\mid A_1A_2\ldots A_{n-1}\right\}
\\
&\geq&
\prod_{i=1}^n\exp\left[-\frac{I_u^{\psi_i^n}(\phi_i)+\frac{\gamma}{n}}{\e}
\right]=\exp\left[-\frac{I_u^{\psi^n}(\phi)+\gamma}{\e} \right].
\end{eqnarray*}

Now for any $r>0$ and $\delta>0$\,,
\begin{eqnarray*}
&&\mathbb{P}\left\{\rho_{0T}(\tilde{u}^{\e,\psi^n}, K_T(r))>\delta \right\}\\
&\leq&\mathbb{P}\Big\{\rho_{0T}(\tilde{u}^{\e,\psi_i^n},
K_{[t_i,t_{i+1}]}(r))>\delta \quad \text{for some } 0\leq
i\leq n-1
\Big\}\\
&\leq& \exp\left[-\frac{r-\gamma}{\e}\right]\,.
\end{eqnarray*}
The proof is complete.
\end{proof}

The lower semicontinuous property  of $I_u^\psi$ for any $\psi\in
C(0, T; H)$ is needed in our approach. We have
\begin{lemma}\label{lem:lower-semicont}
For any $T>0$\,, $I_u^\psi(\phi)$ is lower semicontinuous in both
$\phi$~and~$\psi\in C(0,T; H)$\,.
\end{lemma}
\begin{proof}
The lower semicontinuous property in $\psi$ is followed by a similar
discussion in the proof of second part of Theorem
\ref{thm:ldp-bar-u} by the fact that both~$\bar{f}(\psi)$ and
$\sqrt{B(\psi)}$ are continuous in $\psi$ in space $C(0, T; H)$\,.
And the lower semicontinuous property in $\phi$ is followed by the
same discussion for the usual evolutionary equation~\cite{CeRoc}\,.
\end{proof}

For every $E\subset C(0, T; H)$, denote by~$\operatorname{Int}(E)$
the interior of~$E$ and $\operatorname{Cl}(E)$ the closure of~$E$.
Then by the above result and the lower semicontinuity
of~$\tilde{I}_u^\psi(\phi)$ in~$\phi$ we then have, for any step
function $\psi\in C(0, T; H)$ in the form~(\ref{e:step-function}).
\begin{corollary}\label{cor:ldp-coro}
For any $E\subset C(0,T; H)$\,,
\begin{eqnarray*}
-\inf_{\phi\in \operatorname{Int}(E)}\tilde{I}_u^\psi(\phi)&\leq&
\liminf_{\e\rightarrow 0}\e\ln\mathbb{P}\{u^{\e, \psi}\in E\}
\\&\leq&
\limsup_{\e\rightarrow 0}\e\ln\mathbb{P}\{u^{\e, \psi}\in E\}\leq
-\inf_{\phi\in \operatorname{Cl}(E)}\tilde{I}_u^\psi(\phi)\,.
\end{eqnarray*}
\end{corollary}

Now for our purpose we need the following result
\begin{theorem}
For any $h\in L^2(0, T; H)$, the skeleton equation~(\ref{e:skeleton}) has a unique solution $\phi\in C(0, T; H)$.
\end{theorem}
\begin{proof}
By the assumptions on~$f$ both~$\bar{f}(\phi)$  and~$\sqrt{B(\phi)}$
are Lipschitz continuous. Then the result follows by a standard
discussion on deterministic \textsc{pde}s~\cite{CeRoc}\,.
\end{proof}
Then we define the rate function $I_u(\phi)=I_u^\phi(\phi)$ for any
$\phi\in C(0, T; H)$. Furthermore, for any~$\psi\in C(0, T; H)$,
the following relation between rate functions $I_u^\psi$~and~$I_u$
is needed to derive the \ldp{} for~$\{u^\e\}_\e$.
\begin{lemma}\label{lem:rate function}
Let $\psi^n$ be a family of step functions uniformly converging to
$\phi\in C(0,T; H)$ as $n\rightarrow\infty$\,. Then there is a
family of functions $\phi^n\in C(0, T; H)$ converging to~$\phi$, such
that
\begin{equation*}
\limsup_{n\rightarrow\infty}I_u^{\psi^n}(\phi^n)\leq I_u(\phi)\,.
\end{equation*}
\end{lemma}
\begin{proof}
Suppose $I_u(\phi)<r<\infty$\,. Otherwise the result is clear.

By the definition of~$I_u(\phi)$, there exists a function~$h$ and a
sequence $h^n\in L^2(0, T; H)$ such that
\begin{equation}\label{e:h-phi}
\phi(t)=S(t)u_0+\int_0^tS(t-s)\bar{f}(\phi(s))\,ds
+\int_0^t\sqrt{B(\phi(s))}h(s)\,ds
\end{equation}
and $h^n\rightarrow h$ in $L^2(0, T; H)$ with
\begin{equation*}
\frac12\int_0^T\|h^n(s)\|_0^2\,ds\leq I_u(\phi)+\frac{1}{n}\,.
\end{equation*}
Then define
\begin{equation*}
\phi^n(t)=S(t)u_0+\int_0^tS(t-s)\bar{f}(\psi^n(s))\,ds
+\int_0^t\sqrt{B(\psi^n(s))}h^n(s)\,ds\,.
\end{equation*}
By the Lipschitz property of~$f$, the definition of~$\sqrt{B(\psi)}$
and that $\psi^n\rightarrow \phi$ in space~$C(0, T; H)$ we
have
\begin{equation*}
\phi^n\rightarrow \phi\quad \text{in } C(0, T ;H) \quad
\text{as } n\rightarrow \infty
\end{equation*}
and
\begin{equation*}
\limsup_{n\rightarrow\infty}I_u^{\psi^n}(\phi^n)\leq I_u(\phi)\,.
\end{equation*}

This completes the proof.
\end{proof}

\section{LDP for slow-fast stochastic partial differential equations}
\label{sec:ldp}

Now we show that $\{u^\e\}_\e$ satisfies the \ldp{} with a good rate
function~$I_u$. To do this we prove a special relationship between~$u^\e$ and~$\tilde{u}^{\e,\psi}$ in space~$C(0, T; H)$. The relationship shows that $\tilde{u}^{\e,\psi}$ is comparable with~$u^\e$ near~$\psi$ and implies the \ldp. Here we do not restrict $\psi\in
C(0, T; H)$ to be a step function.

For this let $U^{\e,\psi}=u^\e-\tilde{u}^{\e,\psi}$, then
\begin{equation}\label{e:Ue}
\dot{U}^{\e, \psi}=AU^{\e, \psi}+f(u^\e, v^\e)-f(\psi,
\tilde{v}^{\e,\psi})\,,\quad U^{\e,\psi}(0)=0\,.
\end{equation}
In the mild sense
\begin{eqnarray*}
U^{\e, \psi}(t)&=&\int_0^tS(t-s)\left[ f(u^\e(s),
v^\e(s))-f(\psi(s), \tilde{v}^{\e,\psi}(s))\right]\,ds\\
&=&\int_0^tS(t-s)\left[ f(u^\e(s), v^\e(s))-f(\psi(s),
v^\e(s))\right]\,ds\\&&{}+\int_0^tS(t-s)\left[ f(\psi(s),
v^\e(s))-f(\psi(s), \tilde{v}^{\e,\psi}(s))\right]\,ds \,.
\end{eqnarray*}
By the assumptions on~$f$ and the analysis on fast motion~$\tilde{v}^{\e,\psi}$
\begin{eqnarray*}
\|f(u^\e, v^\e)-f(\psi, v^\e)\|_0&\leq& L_f\|u^\e-\psi\|_0\,,\\
\|f(\psi, v^\e)-f(\psi, \tilde{v}^{\e,\psi})\|_0&\leq&
L_f\|v^\e-\tilde{v}^{\e, \psi}\|_0\leq L_fD_v\|u^\e-\psi\|_0\,.
\end{eqnarray*}
Then there is a positive constant~$C_T$ such that
\begin{equation}\label{e:relation1}
\|u^\e-\tilde{u}^{\e,\psi}\|_{C(0, T; H)}\leq C_T\|u^\e-\psi\|_{C(0,
T; H)}\,.
\end{equation}
Rewrite~(\ref{e:Ue}) as
\begin{eqnarray*}
\dot{U}^{\e,\psi}&=&AU^{\e, \psi}+f(\tilde{u}^{\e,\psi},
\tilde{v}^{\e,\psi})-f(\psi,\tilde{v}^{\e,\psi})+f(u^\e,
v^\e)-f(\tilde{u}^{\e,\psi}, v^\e)
\\&&{}+f(\tilde{u}^{\e,\psi}, v^\e)-f(\tilde{u}^{\e,\psi}, \tilde{v}^{\e,\psi})
\end{eqnarray*}
Then by the Lipschitz property of~$f$ we have one positive constant~$L$ depending on~$L_f$ and~$D_v$ such that
\begin{equation*}
\frac12\frac{d}{dt}\|U^{\e, \psi}\|_0^2\leq
L\|U^{\e,\psi}\|_0^2+L_f(1+D_v)\|\tilde{u}^{\e,\psi}-\psi\|_0^2
\end{equation*}
which yields that
\begin{equation}\label{e:relation2}
\|\tilde{u}^{\e,\psi}-u^\e\|_{C(0, T; H)} \leq
C_T\|\tilde{u}^{\e,\psi}-\psi\|_{C(0,T:H)}
\end{equation}
for some positive constant~$C_T$. Now by relations~(\ref{e:relation1}) and~(\ref{e:relation2}) we prove the main
result.

\begin{proof}[Proof of Theorem~\ref{thm:main}]
We follow Freidlin and Wentzell's~\cite{FW98}
approach to obtain \ldp{} for slow-fast random ordinary differential
equations. For any~$\phi\in C(0, T; H)$ with
$I_u(\phi)<\infty$ and any~$\gamma,\delta>0$\,, by Lemma~\ref{lem:rate function}, we can choose a step function~$\psi^n$ and
a function~$\phi^n$ such that
\begin{equation*}
\rho_{0T}(\phi^n, \phi)<\frac{1}{n}\,,\quad \max_{0\leq t\leq
T}\|\psi^n(t)-\phi(t)\|_0<\frac{1}{n}
\quad\text{and}\quad
\tilde{I}_u^{\psi^n}(\phi^n)<I_u(\phi)+\frac{1}{n}\,.
\end{equation*}
Now by~(\ref{e:relation2})
\begin{eqnarray*}
\rho_{0T}(u^\e,\phi)&\leq&
\rho_{0T}(\tilde{u}^{\e,\psi^n},\phi^n)+\rho_{0T}(u^\e,
\tilde{u}^{\e,\psi^n})+\rho_{0T}(\phi^n,\phi)\\
&\leq&
\rho_{0T}(\tilde{u}^{\e,\psi^n},\phi^n)+C_T\rho_{0T}(\tilde{u}^{\e,\psi^n},
\psi^n)
+\rho_{0T}(\phi^n, \phi) \\
&\leq& \rho_{0T}(\tilde{u}^{\e,
\psi^n},\phi^n)+C_T\rho_{0T}(\tilde{u}^{\e,\psi^n}, \phi^n)
\\&&{}
+C_T\max_{0\leq t\leq T}\|\psi^n(t)-\phi^n(t)\|_0+\rho_{0T}(\phi^n, \phi) \\
 &\leq&
(1+C_T)\rho_{0T}(\tilde{u}^{\e,\psi^n},\phi^n)+(2C_T+1)\frac{1}{n}\,.
\end{eqnarray*}
Then for any $\delta>0$ and $\gamma>0$\,, by Theorem~\ref{thm:LDP-2}, there is an $\e_1>0$ such that for any
$0<\e<\e_1$ we have the following lower bound estimate by
choosing~$n$ large enough
\begin{eqnarray*}
&&\mathbb{P}\left\{\rho_{0T}(u^\e, \phi)\leq\delta\right\}\\
&\geq&\mathbb{P}\left\{(1+C_T)\rho_{0T}(\tilde{u}^{\e,\psi^n},\phi^n)\leq \delta/2\right\}\\
&\geq &\exp\left\{-\frac{I_u^{\psi^n}(\phi^n)+\gamma}{\e}
\right\}\geq \exp\left\{-\frac{I_u(\phi)+2\gamma}{\e} \right\}\,.
\end{eqnarray*}

Now we prove the upper bound estimate. First, by the exponential
tightness of~$\{u^\e\}_\e$ in space~$C(0, T; H)$, for any $r>0$
there is a compact set~$K_r$ such that
\begin{equation*}
\limsup_{\e\rightarrow 0}\e\ln\mathbb{P}\left\{u^\e\in
K^c_r\right\}\leq -r\,.
\end{equation*}
Then for any $\gamma>0$ there is an $\e_2>0$ such that for any
$0<\e<\e_2$
\begin{equation*}
\mathbb{P}\left\{u^\e\in K_r^c\right\}\leq
\exp\left\{-\frac{r-\gamma}{\e} \right\}\,.
\end{equation*}
As $K_r$ is compact, choose a finite $\delta'$-net in~$K_r$ with
$\delta'<\delta$ and let $\phi_1,\phi_2,\ldots, \phi_n$ be the elements of
this net, not belonging to~$K_T(r)$. Then
\begin{equation*}
\mathbb{P}\left\{\rho_{0T}(u^\e, K_T(r))>\delta\right\}\leq
\sum_{i=1}^n\mathbb{P}\left\{\rho_{0T}(u^\e,\phi_i)<\delta'\right\}
+\mathbb{P}\left\{u^\e\in K_r^c\right\}
\end{equation*}
Now we choose step functions $\psi_1,\psi_2,\ldots,\psi_n$ such that
\begin{equation*}
\rho_{0T}(\psi_i, \phi_i)<\delta'\,, \quad i=1,2,\ldots, n\,.
\end{equation*}
Then by the inequality~(\ref{e:relation2}) for $i=1,2,\ldots, n$
\begin{equation*}
\mathbb{P}\left\{\rho_{0T}(u^\e, \phi_i)<\delta'\right\}\leq
\mathbb{P}\Big\{\rho_{0T}(\tilde{u}^{\e, \psi_i},
\phi_i)<2(C_T+1)\delta'\Big\}\,.
\end{equation*}
And by Corollary~\ref{cor:ldp-coro} for $i=1,2,\ldots, n$ we
have
\begin{eqnarray*}
&&\mathbb{P}\Big\{\rho_{0T}(\tilde{u}^{\e, \psi_i},
\phi_i)<2(C_T+1)\delta'\Big\}\\&\leq&
\exp\left\{-\frac{1}{\e}\left[\inf\left\{I_u^{\psi_i}(\phi):
\rho_{0T}(\phi,
\phi_i)<2(C_T+1)\delta'\right\}-\gamma\right]\right\}\,.
\end{eqnarray*}
By the semicontinuity of the functional $I_u^\psi$ in $\psi$, Lemma
\ref{lem:lower-semicont}\,, that for any~$\gamma>0$\,, there
is~$\delta'$ such that
\begin{equation*}
I_u^{\psi_i}(\phi)>r-\gamma/2 \quad \text{for}\quad
\rho_{0T}(\phi, \phi_i)<2(C_T+1)\delta' \quad \text{and} \quad
I^{\phi_i}_u(\phi_i)>r\,.
\end{equation*}
Then  by the choice of $\phi_i\notin K_T(r)$ we have
\begin{equation*}
\mathbb{P}\Big\{\rho_{0T}(\tilde{u}^{\e, \psi_i},
\phi_i)<2(C_T+1)\delta'\Big\}\leq
\exp\left\{-\frac{r-\gamma}{\e}\right\}
\end{equation*}
for $\e$ small enough. The proof is complete.
\end{proof}

\section{An example of slow-fast stochastic reaction-diffusion equation}
\label{sec:esfsrd}

Next we consider the following slow-fast \spde{} on the domain $(0,
L)$ with zero Dirichlet boundary condition
\begin{eqnarray}
&&\p_t u^\e=\p_{xx}u^\e+\lambda\sin u^\e- v^\e\,, \quad u^\e(0)=u_0\,, \label{e:exampleu}\\
&&\e\p_tv^\e=\p_{xx}v^\e-v^\e+u^\e+
   \sqrt{\epsilon}\sigma\p_tW(t)\,, \quad v^\e(0)=v_0
   \label{e:examplev}
 \end{eqnarray}
where $W$ is $L^2(0, L)$-valued $Q$-Wiener process and $\lambda,\sigma>0$ are constants.  As usual, the small parameter~$\epsilon$ measures the separation of time scales between the fast modes~$v$ and the slow modes~$u$.  For small~$\epsilon$ our \ldp\ theory gives the \spde~\eqref{e:example-averaged} as an appropriate (weak) model for the stochastic dynamics of the slow modes~$u$.

Apply the \ldp\ theory to this system.  Note that the nonlinear reaction\slash interaction function $f(u, v)=\lambda\sin u-v$ is Lipschitz continuous. Denote the operator $A=\p_{xx}$ with zero boundary condition on $(0, L)$\,. Now for fixed~$u$ the fast system (\ref{e:examplev}) has a unique stationary solution $\eta^{\e, u}$ with distribution
\begin{equation*}
\mu_u=\mathcal{N}\left((I-A)^{-1}u,
\sigma^2\frac{(I-A)^{-1}Q}{2}\right)\,.
\end{equation*}
Then
\begin{equation*}
\bar{f}(u)=\lambda\sin u-(I-A)^{-1}u\,.
\end{equation*}
Let $\eta^u$ be the stationary solution of
\begin{equation*}
\p_tv=\p_{xx}v-v+u+ \sigma\p_tW(t)
\end{equation*}
for any fixed $u\in L^2(0,L)$\,. Then $\eta^u$ distributes as
$\mu_u$ and
\begin{equation*}
B(u)=2\mathbb{E}\int_0^\infty\big[\eta^u(t)-(I-A)^{-1}u\big]\otimes
\big[\eta^u(0)-(I-A)^{-1}u\big]dt\,.
\end{equation*}
Noticing that
\begin{equation*}
\mathbb{E}\eta^u(t)\otimes\eta^u(0)=\sigma^2\exp\left\{ -(I-A)t
\right\}\left[ \frac{(I-A)^{-1}Q}{2}\right]+(I-A)^{-2}u\otimes u
\end{equation*}
 we then have
\begin{equation*}
\sqrt{B(u)}=(I-A)^{-1}\sigma\sqrt{Q}
\end{equation*}
which is independent of $u$ and satisfies assumption~$\textbf{H}_5$. By Theorem~\ref{thm:main}, $\{u^{\e}\}$~satisfies \ldp{} on space $C(0, T; H)$ with good rate function
\begin{equation*}
I_u(\phi)=\inf_{h\in L^2(0,T;
H)}\left\{\frac12\int_0^T\|h(s)\|_0^2\,ds: \phi=\phi^h\right\}
\end{equation*}
with $\phi^h$ solving
\begin{equation*}
\dot{\phi}=A\phi^h+\lambda\sin\phi-(I-A)^{-1}\phi+(I-A)^{-1}\sigma\sqrt{Q}h\,,\quad
\phi(0)=u_0\,.
\end{equation*}
Furthermore, the rate function is
\begin{equation}\label{e:example-rate}
I_u(\phi)=\frac{1}{2}\int_0^T\left\|\frac{I-A}{\sqrt{Q}}
\left[\dot{\phi}(s)-A\phi(s)-\lambda\sin\phi+(I-A)^{-1}\phi\right]\right\|_0^2\,ds
\end{equation}
for $\phi$ is absolute continuous. Otherwise $I_u(\phi)=\infty$\,.

Now we write out the averaged equation plus the deviation for (\ref{e:exampleu})--(\ref{e:examplev}) as
\begin{equation}\label{e:example-averaged}
d\bar{u}^\e=[A\bar{u}^\e+\lambda\sin \bar{u}^\e-(I-A)^{-1}\bar{u}^\e]\,dt
+\sqrt{\e}\sigma(I-A)^{-1}\sqrt{Q}d\overline{W}(t)
\end{equation}
Then by the \ldp{} for stochastic evolutionary equation~\cite{PZ92},
$\{\bar{u}^\e\}$ satisfies \ldp{} with rate function $I_u(\phi)$
defined in~(\ref{e:example-rate}). This shows that the averaged
equation plus deviation (\ref{e:example-averaged}) does describe the
metastability of $\{u^\e\}_\e$ solving~(\ref{e:exampleu}) for
small~$\e$. Moreover, for large enough parameter~$\lambda$, the
\spde\ model~(\ref{e:example-averaged}) has two stable states near
zero for $\e=0$\,. When $\e\neq 0$\,, noise causes orbits near one
stable state to the position near the other one which shows the
metastability of the system~(\ref{e:example-averaged}). So the
slow-fast stochastic~(\ref{e:exampleu})--(\ref{e:examplev}) also has
such metastability described by system~(\ref{e:example-averaged}).
The description of such tunnelling of the orbit needs detail
analysis by the \ldp{} which is left for future work.

\section{Stochastic centre manifold models confirm the LDP example}
\label{sec:ssmeg}

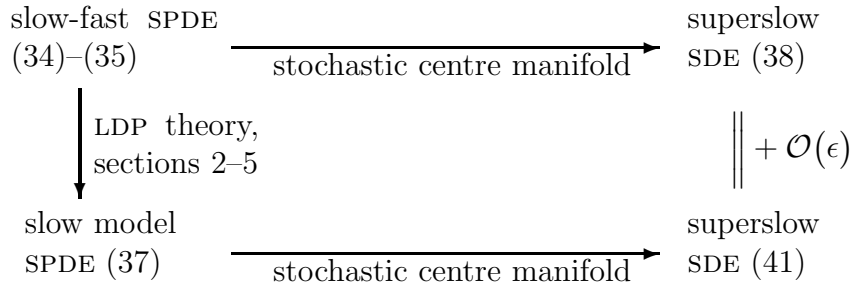
\begin{figure}
\centering
\setlength{\unitlength}{1ex}
\begin{picture}(61,20)
%\put(0,0){\framebox(61,20){}}
\thicklines
\put(0,17){\parbox{15ex}{slow-fast \spde\ \eqref{e:exampleu}--\eqref{e:examplev}}}
\put(1,2){\parbox{11ex}{slow model \spde\ \eqref{e:example-averaged}}}
\put(5,14){\vector(0,-1){8}}
\put(6,9.5){\parbox{12ex}{\ldp\ theory, sections~\ref{sec:pre}--\ref{sec:ldp}}}
\put(16,17){\vector(1,0){31}}
\put(19,15){stochastic centre manifold}
\put(16,2){\vector(1,0){31}}
\put(19,0){stochastic centre manifold}
\put(49,17){\parbox{10ex}{superslow \sde\ \eqref{eq:sfssm}}}
\put(49,2){\parbox{10ex}{superslow \sde\ \eqref{eq:ldpssm}}}
\put(52,9){$\bigg\Vert+\Ord\epsilon$}
\end{picture}
\caption{this schematic diagram shows how the construction of stochastic slow manifold models, in Section~\ref{sec:ssmeg}, verifies by the approximate equality of the two superslow models that the \ldp\ theory correctly captures the slow dynamics of the slow-fast \spde~\eqref{e:exampleu}--\eqref{e:examplev}.}
\label{fig:scheme}
\end{figure}

This section uses the example slow-fast \spde\ \eqref{e:exampleu}--\eqref{e:examplev} to verify the \ldp\ theory in a significant parameter regime.  Without loss of generality set the non-dimensional domain length $L=\pi$.  In the absence of noise the \spde\ \eqref{e:exampleu}--\eqref{e:examplev} undergoes a deterministic pitchfork bifurcation from the trivial field $u=v=0$, as parameter~$\lambda$ crosses the critical value~$3/2$, to two nontrivial fixed points $u\approx \sqrt{\lambda-3/2}\,\sin x$\,.  Consequently, with noise, a stochastic pitchfork bifurcation takes place in the vicinity of parameter $\lambda \approx 3/2$ \cite[e.g.]{Berglund03, Srinamachchivaya90, Blomker06}.  As shown schematically in Figure~\ref{fig:scheme}, here we establish that the stochastic bifurcation dynamics of the original slow-fast \spde\ \eqref{e:exampleu}--\eqref{e:examplev} and that of the \ldp\ averaged \spde~\eqref{e:example-averaged} are identical to the expected \Ord\epsilon~error.

As indicated in  Figure~\ref{fig:scheme}, we make a wide ranging comparison of the dynamics by constructing and comparing the stochastic centre manifolds, and the evolution thereon, of the two \spde\ systems near this stochastic bifurcation.  We explore dynamics in the vicinity of the stochastic bifurcation by setting the parameter $\lambda=\rat32+\lambda'$ for any small enough bifurcation parameter~$\lambda'$.  For small~$\lambda'$ and small~$\epsilon$ there are three time scales in the example fast-slow \spde~\eqref{e:exampleu}--\eqref{e:examplev}: the $v$~modes quasi-equilibrate on the fast time scale of~\Ord\epsilon; almost all of the $u$~modes evolve on the slow time scale of~\Ord1; but the $\sin x$ mode in~$u$ evolves on the superslow, long time scale of~\Ord{1/\lambda'}.  The \ldp\ averaged \spde\ has just the latter two time scales in this parameter regime.  The interactions among \emph{three} time scales is a major complicating factor in constructing the stochastic centre manifold and the evolution thereon. Because of the three time scales, and to be consistent with the terminology of earlier sections, we henceforth refer to the stochastic superslow manifold as it is the evolution on the superslow, long time scales of~\Ord{1/\lambda'} that we encompass and compare in this section.

The stochastic superslow manifold (\ssm) is based from the linear dynamics exactly at critical and, for simplicity, based from no noise~\cite[e.g.]{Roberts06k, Arnold03}. Exactly at critical, and with no noise, $\lambda'=\sigma=0$, both the \spde{}s, \eqref{e:exampleu}--\eqref{e:examplev} and \spde~\eqref{e:example-averaged}, have centre subspaces about the origin: $u=a\sin x$\,, $v=\rat12a\sin x$ and $\bar u=\bar a\sin x$\,, respectively.  Stochastic centre manifold theory~\cite[e.g.]{Boxler89, Arnold03} then asserts that in a domain of small but finite amplitudes $a$~and~$\bar a$, and small but finite noise~$\sigma$ and small but finite parameter~$\lambda'$, there exists an emergent \ssm: solutions are attracted to the \ssm\ roughly as~$\exp(-\rat{27}{10}t)$, and then evolve on the superslow long time scale.  We compare the slow-fast \spde~\eqref{e:exampleu}--\eqref{e:examplev} with the \ldp\ averaged \spde~\eqref{e:example-averaged} via construction of their \ssm\ models.

Constructing \ssm{}s and the evolution thereon has many technical
challenges reported in detail elsewhere~\cite[e.g.]{Chao95,
Roberts06k, Sch08, WangRoberts08}.  The technicalities are even more
challenging here due to the three times scales in the slow-fast
\spde~\eqref{e:exampleu}--\eqref{e:examplev} when, as we assume, the
parameter~$\epsilon$ is small.  The construction procedure used
herein is detailed in a separate technical report~\cite{Roberts09b}
that all can check, reproduce and perhaps modify to other problems
in the same class.  For our purposes we appeal to a little more of
the theory of stochastic centre manifolds:  Arnold~\cite{Arnold03},
building on the work of Boxler~\cite{Boxler91}, assures us that if
the \spde{}s are satisfied to some order of residual in the small
parameters, then the \ssm{}s and the evolution thereon are
constructed to the same order of error.  Thus one may confirm the
veracity of the following \ssm{}s by substituting the expressions
into the \spde{}s and verifying that the residuals are as
asymptotically small as required.

To reduce complicating detail but retain significant information in the example, we truncate the noise to its first three spatial modes: $\dot W=\phi_1\sin x+ \phi_2\sin 2x+ \phi_3\sin 3x$ where $\phi_i$~denote formal derivatives of independent Wiener processes.  Including more noise modes appears to just greatly increase detail, without adding any significant change to the nature of the interactions seen among these three modes.

\paragraph{Slow-fast SPDE \eqref{e:exampleu}--\eqref{e:examplev}}
In six iterations, computer algebra~\cite{Roberts09b} constructs the stochastic superslow manifold model. In terms of the superslow evolving amplitude~$a(t)$, where $u\approx a\sin x$ and $v\approx \rat12 a\sin x$\,,  the stochastic bifurcation \sde\ for the amplitude is
\begin{align}
\dot a={}&\lambda'(1+\rat14\epsilon)a
-(\rat3{16}+\rat18\lambda'+\rat3{64}\epsilon)a^3
+\rat{91}{9728}a^5
\nonumber\\&{}
-\sqrt\epsilon\sigma\left[ (\rat12+\rat18\epsilon
%-\rat14\epsilon\lambda' +\rat9{64}\epsilon a^2
)\phi_1
+\rat3{1216}%+\rat3{4864}\epsilon)
a^2\phi_3 \right]
\nonumber\\&{}
+\epsilon\sigma^2a\left[-\rat1{180}\phi_2\Za2\phi_2
+\rat3{1216}\phi_1\Za3\phi_3
-\rat3{6080}\phi_3\Za3\phi_3
\right]
\nonumber\\&{}
+\Ord{a^6+\lambda'^3+\epsilon^3+\sigma^6}
\label{eq:sfssm}
\end{align}
In this and other expressions, convolutions $\Z{-\alpha}\phi =\int_0^\infty e^{-\alpha s}\phi(t-s)\,ds$\,. The corresponding \ssm\ has slow field
\begin{align}
 u={}& a\sin x
+\rat5{608}a^3\sin 3x
+\rat12\sqrt\epsilon\sigma\sin x\,\Zb1\phi_1
\nonumber\\&{}
-\sqrt\epsilon\sigma\left[
\rat15\sin2x\Big(\Za2-\Zb2 \Big)\phi_2
+\rat1{10}\sin3x \Big(\Za3-\Zb3 \Big)\phi_3 \right]
\nonumber\\&{}
+\Ord{a^4+\lambda'^2+\epsilon^2+\sigma^4},
\label{eq:sfssmu}
\end{align}
and the fast field
\begin{align}
v={}&\rat12a\sin x
+\rat1{1216}a^3\sin3x
\nonumber\\&{}
+\frac\sigma{\sqrt\epsilon}\sin x\left[(1+\rat14\epsilon)\Zb1 +\rat12\Zb1\Zb1\right]\phi_1
\nonumber\\&{}
+\frac\sigma{\sqrt\epsilon}
\sin2x\left[(1+\rat1{25}\epsilon)\Zb2 -\epsilon\rat1{25}\Za2 +\rat15\Zb2\Zb2 \right]\phi_2
\nonumber\\&{}
+\frac\sigma{\sqrt\epsilon}
\sin3x\left[(1+\rat1{100}\epsilon)\Zb3 -\epsilon\rat1{100}\Za3 +\rat1{10}\Zb3\Zb3 \right]\phi_3
\nonumber\\&{}
+\Ord{a^4+\lambda'^2+\epsilon^2+\sigma^4}.
\end{align}
This fast field~$v$ has large $\Ord1$~fluctuations, through terms like $\frac1{\sqrt\epsilon}\Zb1$, because convolution over the fast time scale, $e^{-\beta t/\epsilon}\star$, is~\Ord{\sqrt\epsilon}.   However, note that repeated such convolution, $e^{-\beta t/\epsilon}\star e^{-\beta t/\epsilon}\star$, is~\Ord{\epsilon^{3/2}} \cite[equation~(27)]{Roberts05c}.

\paragraph{LDP averaged SPDE \eqref{e:example-averaged}}
In just four iterations, computer algebra~\cite{Roberts09b} constructs the stochastic superslow evolution to be
\begin{align}
\dot {\bar a}={}&\lambda'{\bar a}-(\rat3{16}+\rat18\lambda'){\bar a}^3
+\rat{91}{9728}{\bar a}^5
\nonumber\\&{}
-\sqrt\epsilon\sigma\left[\rat12\phi_1+\rat3{1216}{\bar a}^2\phi_3\right]
\nonumber\\&{}
+\epsilon\sigma^2\bar a\left[-\rat1{180}\phi_2\Za2\phi_2
+\rat3{1216}\phi_1\Za3\phi_3
-\rat3{6080}\phi_3\Za3\phi_3
\right]
\nonumber\\&{}
+\Ord{{\bar a}^6+\lambda'^3+\epsilon^3+\sigma^6}.
\label{eq:ldpssm}
\end{align}
The corresponding \ssm\ is
\begin{align}
{\bar u}={}&
{\bar a}\sin x +\rat5{608}{\bar a}^3\sin3x
\nonumber\\&{}
-\sqrt\epsilon\sigma\left[
\rat15\sin2x\,\Za2\phi_2 +\rat1{10}\sin3x\,\Za3\phi_3 \right]
%\right.\\&\left.{}
%+\lambda'\left[ \rat15\sin2x\,\Za2\Za2\phi_2 +\rat1{10}\sin3x\,\Za3\Za3\phi_3 \right]
\nonumber\\&{}
+\Ord{{\bar a}^4+\lambda'^2+\epsilon^2+\sigma^4}.
\label{eq:ldpssmu}
\end{align}

\paragraph{Compare the two superslow models}  First compare the slow field~$u$ for the slow-fast \spde,~\eqref{eq:sfssmu}, with the slow field for the \ldp\ averaged \spde,~\eqref{eq:ldpssmu}.  The differences are the $\Ord\epsilon$~terms in the fast time convolutions~$\sqrt\epsilon e^{-\beta t/\epsilon}\star$.  Since $\bar u$ is $u$~averaged over fast fluctuations, these differences are acceptable, and also ensure that the two amplitudes correspond: $\bar a=a+\Ord\epsilon$.

Second, compare the evolution of the amplitudes, \eqref{eq:sfssm}~and~\eqref{eq:ldpssm}.  The only differences are in terms~\Ord\epsilon, as indicated schematically on the right-hand side of Figure~\ref{fig:scheme}. Thus the dynamics of the two superslow \sde\ are within the claimed accuracy of the \ldp.  We conclude that this section verifies that in a parameter regime at least near the stochastic bifurcation, the \ldp\ averaging approximation is correct. Although the \ldp\ averaging only assures us that the slow model is correct in a weak sense, the strong identity between the convolutions appearing in the two \ssm{}s, \eqref{eq:sfssm}--\eqref{eq:sfssmu} and \eqref{eq:ldpssm}--\eqref{eq:ldpssmu}, suggests the correspondence between the \ldp\  averaged system and the original system is generally stronger.

\paragraph{Acknowledgements} This research is supported by the
Australian Research Council grants DP0774311 and DP0988738 and NSFC grant 10701072.

\end{document}